\documentclass[12pt,a4paper,oneside]{amsart}
\usepackage[german]{babel}
\usepackage{amssymb}
\usepackage{amsthm}
\usepackage{booktabs}
\usepackage{float}
\usepackage{hyperref}

\theoremstyle{plain}
\newtheorem{thm}{Theorem}
\newtheorem{cor}{Korollar}
\theoremstyle{definition}
\newtheorem*{defn}{Definition}

\begin{document}
\pagestyle{empty}
\title[]{Lineare Rekurrenzen, Potenzreihen\\ und ihre erzeugenden Funktionen}
\author{Ralf Stephan}
\begin{abstract}
Diese kurze Einf\"uhrung in Theorie und Berechnung linearer Rekurrenzen
versucht, eine L\"ucke in der Literatur zu f\"ullen. Zu diesem Zweck sind
viele ausf\"uhrliche Beispiele angegeben.

This short introduction to theory and usage of linear recurrences tries
to fill a gap in the literature by giving many extensive examples.
\end{abstract}

\maketitle
\parindent0pt
\parskip5pt

\section{Vorwort}
In diesem Artikel m\"ochte ich mich mit linearen Rekursionsgleichungen und den erzeugenden Funktionen ihrer Folgen befassen. Wenngleich die zu verwendenden Techniken nicht neu oder schwierig sind, so bleibt die Behandlung des Themas \"ublicherweise au\ss{}erhalb der Lehrb\"ucher. Das ist auch deswegen bedauerlich, weil sich ganzzahlige Folgen wie die Fibonaccizahlen u.\"a. ausdauernder Beliebtheit erfreuen und das Rad der Berechnung ihrer geschlossenen Form trotzdem st\"andig neu erfunden wird. Dar\"uber hinaus ist das Rechnen mit erzeugenden Funktionen ein unabdingbarer Bestandteil der Kombinatorik.

Im ersten Teil werden wir sehen, dass die durch lineare Rekurrenzen definierten Zahlenfolgen, wenn sie als Koeffizienten in Potenzreihen eingesetzt werden, zu erzeugenden (generierenden) Funktionen f\"uhren, die rational sind. Umgekehrt erzeugt jede rationale Funktion eine Potenzreihe, deren Koeffizienten mindestens einer linearen Rekurrenz gen\"ugen.

Im darauf folgenden Teil wird ersichtlich werden, dass diese Zuordnung mehrere Vorteile hat, insbesondere bei der Berechnung der sogenannten geschlossenen Form der erzeugten Zahlenfolge. Wir hoffen, dass die vermittelten Techniken auch anhand der gro\ss{}en Menge von Beispielen deutlich werden. Obwohl es m\"oglich w\"are, alle Methoden in Software zu realisieren, ist dies unseres Wissens (2007) noch nirgendwo implementiert, daher sollten die Beispiele auch als Testf\"alle f\"ur solche Software geeignet sein.

\medskip
R. Stephan, 2007

\part{Einf\"uhrung und etwas Theorie}
\section{Ein Beispiel: Lukaszahlen}
Dieses Kapitel ist ganz \"ahnlich gehalten wie ein entsprechendes Kapitel von Wilf in \emph{Generatingfunctionology} \"uber die Fibonaccizahlen.

Betrachten wir die Folge der Lucaszahlen. Diese sind, wie die Fibonaccizahlen, definiert durch die lineare Rekurrenz
\begin{equation}\label{rgf-lr}
a_{n+2} = a_{n+1} + a_n
\end{equation}
haben jedoch die Anfangswerte $a_0=2$, $a_1=1$, und ergeben daher eine andere Zahlenfolge:
$$
a_n = {2, 1, 3, 4, 7, 11, 18, 29, 47, 76, 123, \ldots}.
$$

\subsection{Erzeugende Funktion}
Wir wollen nun eine \emph{geschlossene Form} f\"ur die Lucaszahlen finden, das 
hei\ss{}t eine Formel in Abh\"angigkeit von $n$, die nur aus einfachen Funktionen 
besteht, wie Addition, Subtraktion, Multiplikation, Wurzel etc.

Dazu setzen wir die Zahlen der Lucasfolge als Koeffizienten in eine formale 
Potenzreihe (mit formal ist gemeint, dass die Frage der Konvergenz dieser Reihe 
f\"ur uns irrelevant ist). Wir nehmen an, dass der Potenzreihe eine wie auch 
immer geartete \emph{erzeugende Funktion} von $x$ entspricht, und nennen sie $L(x)$:
$$
L(x) = \,\sum_{n=0}^{\infty} a_{n}x^{n} = a_{0}+a_{1}x+a_{2}x^{2}+a_{3}x^{3}+\cdots
= 2 + x + 3x^2 + 4x^3 + 7x^4 + \cdots
$$
Um die geschlossene Form zu finden, ben\"otigen wir erst die erzeugende Funktion 
$L(x)$. Dazu multiplizieren wir die Rekurrenz (\eqref{rgf-lr}) mit $x_n$ und 
summieren ins Unendliche:
$$
\sum_{n=0}^{\infty} a_{n+2}x^{n} = \sum_{n=0}^{\infty} a_{n+1}x^{n} + \sum_{n=0}^{\infty} a_{n}x^{n}
$$
Halt! h\"ore ich die Leser rufen, was ist das? Hierbei handelt es sich um eine ungew\"ohnliche Manipulation, die aber v\"ollig korrekt ist, wenn wir sie in jedem Glied der drei Terme betrachten:
$$
a_{2}+a_{3}x+a_{4}x^{2}+a_{5}x^{3}+\cdots 
= a_{1}+a_{2}x+a_{3}x^{2}+a_{4}x^{3}+\cdots 
+ a_{0}+a_{1}x+a_{2}x^{2}+a_{3}x^{3}+\cdots
$$
Da aber diese Teilfolgen so durch $L(x)$ beschrieben werden k\"onnen: 
$$
a_{0}+a_{1}x+a_{2}x^{2}+a_{3}x^{3}+\cdots
= L(x),
$$
$$
a_{1}+a_{2}x+a_{3}x^{2}+a_{4}x^{3}+\cdots
= \frac{L(x)-a_{0}}{x},
$$
$$
a_{2}+a_{3}x+a_{4}x^{2}+a_{5}x^{3}+\cdots
= \frac{\frac{L(x)-a_{0}}{x}-a_{1}}{x},
$$
bekommt die Rekurrenz \eqref{rgf-lr} die Form
$$
\frac{\frac{L(x)-a_{0}}{x}-a_{1}}{x}
= \frac{L(x)-a_{0}}{x} + L(x),
$$
woraus mit den entsprechenden Startwerten folgt
$$
L(x)=\frac{2-x}{1-x-x^{2}}.
$$

\subsection{Geschlossene Form}
Da sich die erzeugende Funktion $L(x)$ in zwei Partialbr\"uche mit einfacherem Nenner zerlegen l\"a\ss{}t, ist davon auszugehen, dass sich auch die Lucaszahlen als Summe zweier Ausdr\"ucke darstellen lassen. Durch L\"osen der quadratischen Gleichung erhalten wir die Nullstellen des Polynoms im Nenner von $L(x)$, und damit dessen Partialbruchzerlegung:
$$
L(x)=\frac{2-x}{1-x-x^{2}}
= \frac{2-x}{(x-r_{1})(x-r_{2})}=\frac{A}{x-r_{1}}+\frac{B}{x-r_{2}},
$$
$$
r_{1}=\frac{1+\sqrt{5}}{2},\quad r_{2}=\frac{1-\sqrt{5}}{2}.
$$
Nach Berechnung von $A$ und $B$ und weiterer Manipulation folgt
$$
L(x)=\frac{1}{1-r_{1}x}+\frac{1}{1-r_{2}x}.
$$
In diese Form gebracht, l\"a\ss{}t sich die Potenzreihen-Identit\"at
$$
\frac{1}{1-cx}=\sum_{n\ge0}c^{n}x^{n}
$$
anwenden und wir erhalten endlich f\"ur die Lucaszahlen $a_n$
$$
\sum_{n\ge0}a_{n}x^{n}=\sum_{n\ge0}(r_{1})^{n}x^{n}+\sum_{n\ge0}(r_{2})^{n}x^{n}
$$
und daher die \emph{Formel von Binet}
$$
a_{n}=r_{1}^{n}+r_{2}^{n}= \left(\frac{1+\sqrt{5}}{2}\right)^{n}+\left(\frac{1-\sqrt{5}}{2}\right)^{n}.
$$

\section{Potenzreihen-Identit\"aten}
Die einfachste Potenzreihe ist die \emph{geometrische Reihe}:
$$
1 + x + x^{2} + x^{3} + x^{4} + \cdots = \sum_{n=0}^\infty x^n=\frac{1}{1-x}.
$$
Sie erzeugt den konstanten Wert 1. Durch Differentiation und anschlie\ss{}ende 
Multiplikation mit $x$ (man spricht auch von der \emph{Anwendung des $xD$-Operators}) erhalten wir jene Funktion, die die Folge $a_n=n$, also die nat\"urlichen Zahlen, erzeugt.
$$
0 + x + 2x^{2} + 3x^{3} + 4x^{4} + \cdots = \sum_{n=0}^\infty nx^n=\frac{x}{(1-x)^{2}}.
$$
Jede weitere Anwendung des $xD$-Operators erzeugt die n\"achste Potenz von $n$.
$$
0 + x + 4x^{2} + 9x^{3} + 16x^{4} + \cdots = \sum_{n=0}^\infty n^2x^n=\frac{x(x+1)}{(1-x)^{3}}.
$$
$$
\sum_{n=0}^\infty n^3x^n=\frac{x(x^{2}+4x+1)}{(1-x)^{4}}.
$$
$$
\sum_{n=0}^\infty n^4x^n=\frac{x(x+1)(x^{2}+10x+1)}{(1-x)^{5}}.
$$
\begin{equation}\label{rgf-p1}
\sum_{n=0}^\infty n^mx^n=(xD)^{m}\frac{1}{(1-x)}.
\end{equation}
Hier bezeichnet $(xD)^m$ die $m$-malige Anwendung des $xD$-Operators. 

Geht es darum, Potenzen einer konstanten Zahl $c$ zu erzeugen, gen\"ugt es, $x$ 
in der geometrischen Reihe durch $cx$ zu ersetzen:
$$
1 + cx + c^{2}x^{2} + c^{3}x^{3} + c^{4}x^{4} + \cdots = \sum_{n=0}^\infty c^nx^n=\frac{1}{1-cx}.
$$
Ebenso gelten die Regeln des $xD$-Operators f\"ur diese Funktion, und wir erhalten
\begin{equation}\label{rgf-p2}
\sum_{n=0}^\infty c^nn^mx^n=(xD)^{m}\frac{1}{(1-cx)}.
\end{equation}

Wir werden sehen, dass alle linearen Rekurrenzen bzw.~ihre rationalen 
erzeugenden Funktionen Zahlenfolgen erzeugen, die sich als Summe von den in 
\eqref{rgf-p1} und \eqref{rgf-p2} erzeugten Ausdr\"ucken darstellen lassen. 
Die Behandlung weiterer Identit\"aten von Potenzreihen ist hier daher gar nicht notwendig.

\section{Homogene und inhomogene Rekurrenzen}
Bevor wir zu den zentralen theoretischen Erkenntnissen kommen, die diesem Buch zugrunde liegen, ist es noch notwendig, zwei Arten von linearen Rekurrenzen zu betrachten. Wir werden sehen, dass nur sogenannte \emph{homogene} Rekurrenzen 
auf einfache Weise handhabbar sind, und dass \emph{inhomogene} Rekurrenzen in diese \"uberf\"uhrt werden k\"onnen. Alle anderen Arten von Rekurrenzen k\"onnen mit den dargestellten Methoden nicht gel\"ost werden. Zun\"achst jedoch:
\begin{defn}
Eine \emph{lineare Rekurrenz} (auch \emph{Rekursion} oder \emph{Differenzengleichung}) besteht aus einer Gleichung der Form
$$
c_ka_{n+k} = c_{k-1}a_{n+k-1} + c_{k-2}a_{n+k-2} + \cdots + c_0a_n + f(n),
$$
und vorgegebenen Anfangswerten 
$$
a_0, a_1, \ldots, a_{k-1}
$$
wobei $n$ eine ganzzahlige Unbekannte ist, $k$ als positive ganze Zahl die \emph{Ordnung} der Rekurrenz bezeichnet, und die $c_i$ vorgegebene ganze Zahlen 
sind. Auch ist die Funktion von $n$, $f(n)$, von einer Form, die sich wiederum 
von einer linearen Rekurrenz darstellen l\"a\ss{}t. Da wir uns nur mit ganzzahligen Folgen befassen, ist $c_k$ gleich Eins.

Eine lineare Rekurrenz ist \emph{homogen}, wenn $f(n)$ gleich Null ist.
\end{defn}

\subsubsection{Beispiele.}
\begin{itemize}
\item Quadratische Rekurrenz: $a_{n+2}=a_{n+1}a_{n},\qquad a_0=2, a_1=3$
\item Nichtlineare Rekurrenz: $a_{n+2}=na_{n+1}-a_{n},\qquad a_0=1, a_1=2$
\item Inhomogene lineare Rekurrenz 2.Ordnung: $\,a_{n+2}=5a_{n+1}-a_{n}-n, \qquad a_0=0, a_1=1$
\item Inhomogene lineare Rekurrenz 3.Ordnung: $a_{n+3}=a_{n+2}+a_{n+1}+a_{n}+1, \qquad a_0=0, a_1=1, a_2=-1$
\item Homogene lineare Rekurrenz 4.Ordnung: $a_{n+4}=2a_{n+3}-a_{n+2}+5a_{n+1}-a_{n}, \qquad a_0=-5, a_1=-1, a_2=3, a_3=-7$
\end{itemize}

\section{Ein Hauptsatz}
F\"ur das Verst\"andnis linearer Rekurrenzen ist es zentral, den folgenden Hauptsatz und seine Folgerungen zu kennen.
\begin{thm}
Jede homogene lineare Rekurrenz erzeugt eine Zahlenfolge, die als Koeffizienten in einer Potenzreihe von einer rationalen Funktion erzeugt wird, wobei der Nenner der erzeugenden Funktion ein Polynom ist, dessen Koeffizienten mit den Koeffizienten der entsprechenden Rekurrenz identisch sind, wenn die Gleichung gleich Null gesetzt wird.
\end{thm}
Damit besteht eine bijektive Zuordnung zwischen rationalen Funktionen und homogenen linearen Rekurrenzen. Als Beispiel siehe den ersten Abschnitt
\"uber die Lucaszahlen. Setzen wir die Differenzengleichung der Lucaszahlen auf Null, lautet sie:
$$
a_{n+2}-a_{n+1}-a_{n}=0
$$
mit den Koeffizienten $1, -1, -1$ gleich denen des Nennerpolynoms von $L(x)$.

Der allgemeine Beweis ist einfach, man verfahre so wie im Beispiel bei der Herleitung der erzeugenden Funktion. Die nachstehenden Folges\"atze ergeben sich unmittelbar:
\begin{cor}
Der Grad des Nennerpolynoms der erzeugenden Funktion ist gleich der Ordnung der entsprechenden linearen Rekurrenz.
\end{cor}
Die wichtigste Schlu\ss{}folgerung ergibt sich jedoch durch Anwendung des 
Fundamentalsatzes der Algebra und der zuvor gefundenen Potenzreihen-Identit\"aten.
\begin{thm}
Jede lineare Rekurrenz hat eine geschlossene Form der Art
$$
a_n=\sum_{j=1}^{D}\frac{1}{z_j^n}\sum_{k=1}^{M_j}c_{j,k}n^{k-1},$$
wobei $D$ die Ordnung der Rekurrenz, $z_j$ eine (auch komplexe) Nullstelle des 
Nennerpolynoms der erzeugenden Funktion, $M_j$ die Vielfachheit dieser 
Nullstelle, und die $c_{j,k}$ rationale Konstanten sind, 
die es schlu\ss{}endlich herauszufinden gilt.
\end{thm}
Siehe den Buchteil ''Ausf\"uhrliche Beispiele'', wo das Finden der geschlossenen Form ein zentrales Thema ist.
\begin{cor}
Das asymptotische Verhalten einer linearen Rekurrenz ist exponentiell. Die Basis der Potenz ist das betragsm\"a\ss{}ig gr\"o\ss{}te Inverse aller Nullstellen des Nennerpolynoms der erzeugenden Funktion.
\end{cor}
Auch dies ergibt sich unmittelbar, und damit kann das asymptotische Verhalten einer linearen Rekurrenz direkt aus der Differenzengleichung berechnet werden, wenn sie in der homogenen Form ist.

Als Beispiel die Formel $2^n+3^n$ mit der Rekurrenz 
$$
a_{n+2}=5a_{n+1}-6a_n,\qquad a_0=2, a_1=5,
$$
der das Nennerpolynom $1-5x+6x^2$ mit den Nullstellen $\frac12$ und $\frac13$ entspricht. Zu beachten ist, dass Nullstellen auch komplex sein k\"onnen.

\part{Ausf\"uhrliche Beispiele}
\section{Die Perrin-Folge}
\"Uber die vor allem zahlentheoretischen Eigenschaften der Perrin-Folge bietet 
die Wikipedia einen Artikel. Wir wollen die geschlossene Form dieser Zahlenfolge ausgehend von ihrer Differenzengleichung
$$
a_{n+3}=a_{n+1}+a_n, \qquad\text{mit $a_0=3, a_1=0, a_2=2$}
$$
herleiten. Sei $P(x)$ die entsprechende erzeugende Funktion
$$
P(x) = \sum_{n\ge0}a_nx^n = 3 + 2x^2 + 3x^3 + 2x^4 + 5x^5 + \cdots,
$$
dann gilt
$$
\frac{P(x) - 3 - 2x^2}{x^3} = \frac{P(x)-3}{x} + P(x),
$$
$$
P(x)-3-2x^2=x^2P(x)-3x^2+x^3P(x),
$$
$$
P(x)=\frac{3-x^2}{1-x^2-x^3}.
$$
Seien $z_1,z_2,z_3$ die L\"osungen der Gleichung $1-x^2-x^3$, dann ist
\begin{equation}\label{rgf-pe1}
P(x)=\frac{3-x^2}{1-x^2-x^3}=\frac{A}{x-z_1}+\frac{B}{x-z_2}+\frac{C}{x-z_3}.
\end{equation}
Zun\"achst wollen wir eine asymptotische Absch\"atzung machen, weiter unten in diesem Kapitel werden wir aber eine M\"oglichkeit zur Bestimmung 
der~$A,B,C$ kennenlernen.

\subsection{Das asymptotische Verhalten der Perrin-Folge}
Die Funktion $P(x)$ ist, wie in Gleichung~\eqref{rgf-pe1} gesehen, die Summe 
dreier Br\"uche, deren Nenner $x-z$ durch Division zu $1-\frac{x}{z}$ umgeformt 
werden kann. Der dabei entstehende  Faktor~$\frac1z$ von $x$ entspricht der 
Konstanten $c$ in einer unserer Potenzreihen-Identit\"aten, und daher k\"onnen wir 
die Perrin-Zahlen als Summe dreier Potenzen, multipliziert mit jeweils einem unbekannten Faktor, darstellen:
\begin{equation}\label{rgf-pe2}
a_n=D(1/z_1)^n+E(1/z_2)^n+F(1/z_3)^n.
\end{equation}
Die drei Nullstellen $z_1,z_2,z_3$ des Nennerpolynoms von $P(x)$ k\"onnen mit 
Taschenrechner oder PC-Software (z.B. ein Computer-Algebra-System) numerisch erhalten werden:
\begin{align*}
z_1 &= 0,754877666\ldots\\
z_2 &= -0,877438833\ldots - 0,744861767\ldots i\\
z_3 &= -0,877438833\ldots + 0,744861767\ldots i
\end{align*}
Egal, welche Faktoren $D,E,F$ gelten, mit zunehmendem $n$ wird sich derjenige 
Term durchsetzen, dessen Basis $\frac1z$ den gr\"o\ss{}ten Betrag besitzt, und wird
das asymptotische Verhalten bestimmen. Der Rechner liefert
\begin{align*}
\Bigl|\frac{1}{z_1}\Bigr| &= 1,324717917\ldots\\
\Bigl|\frac{1}{z_2}\Bigr| &= 0,868936862\ldots\\
\Bigl|\frac{1}{z_3}\Bigr| &= 0,868936862\ldots
\end{align*}
und daher
\begin{equation}\label{rgf-pe3}
a_n \sim D\cdot1.324717917^n.
\end{equation}

\subsection{Partialbruchzerlegung durch Koeffizientenvergleich}
Die Berechnung der Faktoren $A,B,C$ in Gleichung~\eqref{rgf-pe1} ist durch 
folgendes Standardverfahren m\"oglich. Die Gleichung wird mit dem Nennerpolynom 
von~$P(x)$ multipliziert.
$$
3-x^2=A(x-z_2)(x-z_3)+B(x-z_1)(x-z_3)+C(x-z_1)(x-z_2)
$$
$$
3-x^2=A(x^2-z_2x-z_3x+z_2z_3)+B(x^2-z_1x-z_3x+z_1z_3)+C(x^2-z_1x-z_2x+z_1z_2)
$$
Dies f\"uhrt bei Ber\"ucksichtigung der jeweiligen Faktoren von 1, $x$ und $x^2$ 
zu einem Gleichungssystem mit drei Unbekannten
\begin{align*}
3 &= z_2z_3A + z_1z_3B + z_1z_2C\\
0 &= (z_2+z_3)A + (z_1+z_3)B + (z_1+z_2)C\\
-1 &= A + B + C
\end{align*}
mit der L\"osung $A=z_1,B=z_2,C=z_3$, woraus wiederum $D=E=F=1$ folgt.
Durch Einsetzen in die Gleichungen~\eqref{rgf-pe2} und~\eqref{rgf-pe3} 
erhalten wir
$$
a_n=(1/z_1)^n+(1/z_2)^n+(1/z_3)^n.
$$
Unsere vorige Absch\"atzung~\eqref{rgf-pe3} ist so gut, dass bei ihr ab der zehnten Perrinzahl nur noch gerundet werden muss:

\floatstyle{boxed}
\restylefloat{table}
\begin{table}[t]
\begin{center}
\begin{tabular}{rrr}
$n$ & $a_n$ & $\left(\frac{1}{z_1}\right)^n$\\
\toprule\\
 0 &  3  &  3.000000\\
 1 &  0  &  1.324717\\
 2 &  2  &  1.754877\\
 3 &  3  &  2.324717\\
 4 &  2  &  3.079595\\
 5 &  5  &  4.079595\\
 6 &  5  &  5.404313\\
 7 &  7  &  7.159191\\
 8 & 10  &  9.483909\\
 9 & 12  & 12.563504\\
 10& 17  & 16.643100\\
 11& 22  & 22.047414\\
 12& 29  & 29.206605\\
 13& 39  & 38.690514\\
 14& 51  & 51.254019\\
 15& 68  & 67.897119\\
 16& 90  & 89.944533\\
 17&119  &119.15113\\
 18&158  &157.84165\\
 19&209  &209.09567\\
 20&277  &276.99279
\end{tabular}
\end{center}
\caption{Wertetabelle der Perrinfolge und ihrer N\"aherung.}
\end{table}

\section{Eine Fibonacci-Teilfolge}
Bei genauerer Betrachtung der Fibonacci-Folge
$$
f_n = {0, 1, 1, 2, 3, 5, 8, 13, 21, 34, 55, 89, 144, ...}
$$
f\"allt auf, dass die Glieder der Teilfolge 
$$
f_{3n+1} = {1, 3, 13, 55, 233, 987, 4181, ...}
$$
alle ungerade zu sein scheinen, was wir zun\"achst ohne Beweis voraussetzen. Angenommen, uns interessiert die Folge, die sich durch Verringerung um~1
und anschlie\ss{}ender Halbierung dieser Werte ergibt, also
$$
a_n=\frac{f_{3n+1}-1}{2} = \{0, 1, 6, 27, 116, 493, 2090, 8855, 37512, 158905, \ldots \},
$$
und wir wollen eine geschlossene Form daf\"ur finden, dann brauchen wir zuerst 
eine Rekurrenz f\"ur $a_n$. Dazu ist es notwendig, auch $a_{n-1}$ und 
$a_{n-2}$ durch $f_n$ auszudr\"ucken:
$$
a_n=\frac{f_{3n+1}-1}{2} \quad\Longrightarrow\quad a_{n-1}=\frac{f_{3n-2}-1}{2},\quad a_{n-2}=\frac{f_{3n-5}-1}{2}.
$$
Durch einfache Manipulationen und Kenntnis der Rekurrenz f\"ur $f_n$ ist es
m\"oglich, $f_{3n+1}$ durch $f_{3n-2}$ und $f_{3n-5}$ auszudr\"ucken:
$$
f_n=f_{n-1}+f_{n-2}
$$
\begin{align*}
f_{3n+1} &= f_{3n}+f_{3n-1}\\
 &= 2f_{3n-1}+f_{3n-2}\\
 &= 3f_{3n-2}+2f_{3n-3}\\
 &= 4f_{3n-2}+f_{3n-3}-f_{3n-4}\\
 &= 4f_{3n-2}+f_{3n-5}
\end{align*}

Daraus folgt wiederum f\"ur $a_n$
$$
2a_n+1 = 4(2a_{n-1}+1) + (2a_{n-2}+1),
$$
und wir erhalten die Rekurrenz
$$
a_n=4a_{n-1}+a_{n-2}+2,\quad a_0=0,a_1=1.
$$
Die Inhomogenit\"at der Rekurrenz spielt keine Rolle bei der Berechnung der 
erzeugenden Funktion $A(x)$, solange sich zus\"atzliche Ausdr\"ucke als Potenzreihe 
darstellen lassen, und wir erhalten
$$
A(x)=\sum_{n\ge0}a_nx^n \quad=\quad 4\frac{A(x)}{x} + \frac{A(x)-x}{x^2} + 2\frac{1}{1-x},
$$
$$
A(x)=\frac{x+x^2}{(1-x)(1-4x-x^2)}.
$$
L\"osung der quadratischen Nennerfaktor-Gleichung liefert den Ansatz f\"ur die Partialbruchzerlegung
$$
A(x)=\frac{x+x^2}{(1-x)(1-4x-x^2)} = \frac{C}{x-r_+} + \frac{D}{x-r_-} + \frac{E}{1-x},
\qquad r_{\pm}=-2\pm\sqrt5.
$$
Mit der Standard-L\"osungsmethode kommen wir \"uber
$$
x+x^2=C(x-r_-)(1-x)+D(x-r_+)(1-x)+E(1-4x-x^2)
$$
$$
x+x^2=C((1+r_-)x-r_--x^2)+D((1+r_+)x-r_+-x^2)+E(1-4x-x^2)
$$
auf das Gleichungssystem
\begin{align*}
 0 &= -r_-C-r_+D+E\\
 1 &= (1+r_-)C+(1+r_+)D-4E\\
 1 &= -C-D-E
\end{align*}
mit der L\"osung
$$
C=\frac{3\sqrt5-5}{20},\quad D=\frac{-3\sqrt5-5}{20},\quad E=-\frac{1}{2},
$$
woraus sich nach Umwandlung der Nenner die gesuchte Form ergibt
$$
a_n=-\tfrac{1}{2}+\tfrac{1}{20}\big((\sqrt5+5)(\sqrt5+2)^n+(-\sqrt5+5)(-\sqrt5+2)^n\big),
$$
mit der Absch\"atzung
$$
a_n\sim-\tfrac{1}{2}+\tfrac{1}{20}(\sqrt5+5)(\sqrt5+2)^n.
$$

\section{Die Partialsummen der Kubikzahlen}
Ein h\"aufig auftretendes Problem sind \emph{Partialsummen}-Folgen und \emph{Differenzen}-Folgen, die sich aus vorgegebenen Zahlenfolgen ableiten. 
Als Operation auf Potenzreihen betrachtet, erh\"alt man die \emph{erste 
Differenz} durch Multiplikation der Erzeugenden mit $(1-x)$:
\begin{align*}
 A(x) &= a_0+a_1x+a_2x^2+a_3x^3+a_4x^4+\cdots\\
 (1-x)A(x) &= a_0+(a_1-a_0)x+(a_2-a_1)x^2+(a_3-a_2)x^3+(a_4-a_3)x^4+\cdots
\end{align*}
Daraus folgt unmittelbar, dass die Funktion $A(x)/(1-x)$ die \emph{erste Partialsumme} der Folge $a_n$ erzeugt:
$$
\frac{A(x)}{1-x} = \sum_{n\ge0}\left(\,\sum_{k=0}^na_k\right)x^n.
$$
Wir wollen dies am Beispiel der Kubikzahlen $0, 1, 8, 27, 64, 125\ldots$
veranschaulichen.  Die Partialsummenfolge der Kubikzahlen
$$
a_n=\sum_{k=0}^nk^3=\{0,1,9,36,100,225\ldots\}
$$
besitzt die erzeugende Funktion (siehe das Kapitel \"uber Potenzreihen-Identit\"aten)
$$
A(x)=\sum_{n\ge0}a_nx^n=\frac{1}{1-x}\cdot\frac{x(x^2+4x+1)}{(1-x)^4}=\frac{x(x^2+4x+1)}{(1-x)^5}.
$$
Die Nullstelle des Nennerpolynoms hat den Wert~1 und ist f\"unffach, das hei\ss{}t, 
in der geschlossenen Form von $a_n$ k\"onnen Potenzen von~$n$ bis~$n^4$ 
auftreten. Der Ansatz f\"ur die Partialbruchzerlegung von~$A(x)$ lautet
\begin{multline*}
\frac{x(x^2+4x+1)}{(1-x)^4}=\\
\frac{A}{1-x}+\frac{Bx}{(1-x)^2}+\frac{Cx(x+1)}{(1-x)^3}+\frac{Dx(x^2+4x+1)}{(1-x)^4}+\frac{Ex(x+1)(x^2+10x+1)}{(1-x)^5},
\end{multline*}
\begin{align*}
 x^3+4x^2+x= &A(1-x)^4+B(-x^4+3x^3-3x^2+x)+C(x^4-x^3-x^2+x)\\
 &+D(-x^4-3x^3+3x^2+x)+E(x^4+11x^3+11x^2+x)
\end{align*}
mit dem Gleichungssystem
\begin{align*}
 0 &= A\\
 1 &= -4A+B+C+D+E\\
 4 &= 6A-3B-C+3D+11E\\
 1 &= -4A+3B-C-3D+11E\\
 0 &= A-B+C-D+E
 \end{align*}
Es hat die L\"osung $C=\frac14,D=\frac12,E=\frac14$ und daher ist
$$
a_n=\sum_{k=0}^nk^3=\frac{1}{4}n^2+\frac{1}{2}n^3+\frac{1}{4}n^4=\frac{n^2(n+1)^2}{4}.
$$
Was auch gleich dem Quadrat der Dreieckszahlen ist.

Diese Identit\"at l\"a\ss{}t sich selbstverst\"andlich auch auf anderem Weg (und in diesem Fall vielleicht sogar schneller) beweisen. Es ging in diesem Kapitel jedoch um die Demonstration der allgemeinen Methode zur Behandlung von Partialsummen, und um mehrfache Nullstellen bei der Partialbruchzerlegung.

\section{Eine vorgegebene Folge}
Es sei das Problem vorgegeben, die geschlossene Form der Folge zu finden, die die erzeugende Funktion
$$
A(x)=\frac{1+2x-x^2}{(1-x)^4(1+x)^2}=1+4x+8x^2+16x^2+25x^3+40x^4+56x^5+\cdots
$$
besitzt (diese Zahlenfolge entsteht bei einer gewissen Zerteilung von Polygonen). Wir k\"onnen diesmal nicht von der Rekurrenz ausgehen, sehen aber sofort anhand der Erzeugenden, dass sie von 6. Ordnung ist, au\ss{}erdem kommt in 
der geschlossenen Form~$n$ bis zur dritten Potenz, sowie Potenzen von $(-1)$ 
multipliziert mit~$n$ bis zur ersten Potenz, vor. Der Ansatz setzt sich daher 
aus den vier Termen~$1$, $n$, $n^2$ und den zwei Termen~$(-1)^n$ und~$n(-1)^n$
als Einzelfunktionen der Partialbruchzerlegung zusammen:
\begin{multline*}
\frac{1+2x-x^2}{(1-x)^4(1+x)^2}=\\
\frac{Cx(x^2+4x+1)}{(1-x)^4}+\frac{Dx(x+1)}{(1-x)^3}+\frac{Ex}{(1-x)^2}+\frac{F}{1-x}+\frac{-Gx}{(1+x)^2}+\frac{H}{1+x}
\end{multline*}
(Benutzen Sie Gleichung \eqref{rgf-p2} mit $c=-1$, um die letzten beiden Terme zu erhalten!)

Der Ansatz f\"uhrt zu dem Gleichungssystem
$$
\begin{pmatrix}0\\0\\0\\-1\\2\\1\end{pmatrix}\quad=\quad
\begin{pmatrix}C&D&E&F&G&H\end{pmatrix}
\begin{pmatrix}
1&-1&1&-1&-1&1\\
6&-2&0&1&4&-3\\
10&0&-2&2&-6&2\\
6&2&0&-2&4&2\\
1&1&1&-1&-1&-3\\
0&0&0&1&0&1
\end{pmatrix}
$$
mit der L\"osung 
$$
\{C=\tfrac1{12},D=\tfrac{3}{4},E=\tfrac{43}{24},F=\tfrac{9}{8},
G=-\tfrac18,H=-\tfrac18\}
$$
und damit der Formel
$$
a_n=\tfrac1{24}\left(2n^3+18n^2+43n+27-3(n+1)(-1)^n\right).
$$

\section{Ein gemeinsamer Teiler}
Die Zahlenfolge mit der Formel $a_n=\text{ggT}\,(n,4)$ ist 4-periodisch:
$$
a_n= \{4, 1, 2, 1, 4, 1, 2, 1, 4, 1, 2, 1, ...\}
$$
Wenngleich die Folge damit eine eindeutige und leicht berechenbare Form hat, w\"urde uns eine geschlossene Form interessieren, die nur aus Potenzen besteht. Dazu versuchen wir zun\"achst, eine Rekurrenz zu finden.

Die Periodizit\"at liefert sofort:
$$
a_n=a_{n-4} \qquad a_0=4,a_1=1,a_2=2,a_3=1
$$
und daher
$$
A(x)=\sum_{n\ge0}a_nx^n=\frac{4+x+2x^2+x^3}{1-x^4}=\frac{C}{x-1}+\frac{D}{x+1}+\frac{E}{x-i}+\frac{F}{x+i}.
$$
Der Ansatz
\begin{multline*}
4+x+2x^2+x^3=\\
C(1+x+x^2+x^3)+D(x^3+x-x^2-1)+E(x^3-x+ix^2-i)+F(x^3-x-ix^2+i)
\end{multline*}
f\"uhrt auf das Gleichungssystem
\begin{align*}
 4 &= C-D-iE+iF\\
 1 &= C+D-E-F\\
 2 &= C-D+iE-iF\\
 1 &= C+D+E+F
\end{align*}
mit der L\"osung $C=2,D=-1,E=\frac{i}2,F=-\frac{i}2$. Wir erhalten die Formel
$$
a_n=\mathbf{ggT}(n,4)\quad=\quad2+(-1)^n+\frac{i^n+(-i)^n}{2}.
$$

\section{Eine Differenzengleichung f\"ur $\cosh$}
Die in den bisherigen Beispielen angewandten Prinzipien f\"uhrten zu Methoden, die zum Ziel hatten, eine bestimmte Art der geschlossenen Form zu erhalten, n\"amlich eine Summe von Potenzen. Umgekehrt ist es kein Problem, von einer solchen Form auszugehen, und zu einer formalen Potenzreihe zu gelangen, von der sich dann eine Rekurrenz ablesen l\"a\ss{}t.

Beispielsweise lautet die Definition f\"ur den \emph{Kosinus Hyperbolicus}
$$
\cosh(z) = \tfrac12e^z+\tfrac12e^{-z}, \qquad z\in\mathbb{C}.
$$
Daraus folgt f\"ur die formal erzeugende Funktion, wir nennen sie $C(x)$,
$$
C(x)=\frac1{2(1-ex)}+\frac1{2(1-x/e)}=
\frac{2-(e+\frac1e)x}{2(1-(e+\tfrac1e)x+x^2)}.
$$
Der Nenner wiederum zeigt den Weg zur Differenzengleichung
$$
\cosh(n+2)=(e+\tfrac1e)\cosh(n+1)-\cosh n.
$$

\section{Alternierende Folgen}
Mit den Methoden des letzten Beispiels k\"onnen auch komplexere Aufgabenstellungen bew\"altigt werden. Sei die Definition gegeben
$$
f(n) = 
\begin{cases} 
 f_1(n), & n\mbox{ gerade;}\\
 f_2(n), & n\mbox{ ungerade.}
\end{cases}
$$
Solche periodisch wechselnden Folgen werden mit dem Ansatz
$$
f(n)=\frac{f_1(n)+f_2(n)}2+\frac{f_1(n)-f_2(n)}2(-1)^n.
$$
bearbeitet. Das Prinzip: ist $n$ gerade, bleibt das zwischen den Br\"uchen 
stehende Vorzeichen ein Plus, und die $f_2$-Terme heben sich gegenseitig auf. 
Ist $n$ ungerade, wird das Vorzeichen ein Minus und die $f_1$-Terme verschwinden.

Ein Beispiel, von der Definition \"uber die Erzeugende zur Rekurrenz:
\begin{align*}
 f_n &=
 \begin{cases}
  2^n,&n\text{ gerade;}\\
  n+1,&n\text{ ungerade.}
 \end{cases}
 \qquad =\tfrac{2^n+n+1}2+\tfrac{2^n-n-1}2(-1)^n\\
 F(x) &=\frac12\left(\frac1{1-2x}+\frac{x}{(1-x)^2}+\frac1{1-x}
+\frac1{1+2x}-\frac{-x}{(1+x)^2}-\frac1{1+x}\right)\\
 &=\frac{-x^4+8x^3+2x^2-2x-1}{4x^6-9x^4+6x^2-1}.
\end{align*}

Das Resultat
$$
f_{n+6}=6f_{n+4}-9f_{n+2}+4f_n
$$
besteht nur aus Folgengliedern, die einer einzigen der in der Definition 
angegebenen Teilfolgen angeh\"oren, gilt aber f\"ur die ganze Folge. Das bedeutet, 
dass durch Halbierung der Indices eine Rekurrenz entsteht, die f\"ur jede der 
beiden Einzelfolgen $2^{2n}$ und $2n$ gilt, und wir eine Methode entdeckt haben, mit der sich eine solche berechnen l\"a\ss{}t.

\section{Weiterf\"uhrende Literatur}
\"Uber den kombinatorischen Aspekt generierender Funktionen siehe man
Wilf, \emph{Generatingfunctionology} und \"uber ihre asymptotische Entwicklung
und vieles mehr sei Graham, Knuth, Ptashnik, \emph{Concrete Mathematics} empfohlen.

\end{document}